\documentclass{opt2020} 

\usepackage{booktabs}
\usepackage{color}

\usepackage{algorithm}
\usepackage{algorithmic}

\usepackage{url}

\usepackage{amssymb}
\usepackage{amsmath}
\usepackage{mathtools}

\newcommand{\trace}{{\rm trace}}

\newcommand{\hess}{\mathrm{Hess}}

\def\grad{\mathop{\rm grad}\nolimits}

\def\bzero{{\mathbf 0}}

\def\bA{\mathbf A}
\def\bB{\mathbf B}

\def\bI{{\mathbf I}}

\def\bU{{\mathbf U}}

\def\bX{{\mathbf X}}
\def\bY{{\mathbf Y}}
\def\bZ{{\mathbf Z}}
\def\bM{{\mathbf M}}

\def\argmin{\mathop{\rm arg\,min}\limits}

\renewcommand{\trace}{{\rm trace}}

\def\min{\mathop{\rm min}\nolimits}

\begin{document}
\title[Riemannian optimization on the simplex of positive definite matrices]{Riemannian optimization on the \\simplex of positive definite matrices
}



\optauthor{\Name{Bamdev Mishra} \Email{bamdevm@microsoft.com}\\
	\addr Microsoft, India\\
	\Name{Hiroyuki Kasai} \Email{hiroyuki.kasai@waseda.jp}\\
	\addr Waseda University, Japan\\
	\Name{Pratik Jawanpuria} \Email{pratik.jawanpuria@microsoft.com}\\
	\addr Microsoft, India
}


\maketitle

\begin{abstract}
In this work, we generalize the probability simplex constraint to matrices, i.e., $\bX_1 + \bX_2 + \ldots + \bX_K = \bI$, where $\bX_i \succeq 0$ is a symmetric positive semidefinite matrix of size $n\times n$ for all $i = \{1,\ldots,K \}$. By assuming positive definiteness of the matrices, we show that the constraint set arising from the matrix simplex has the structure of a smooth Riemannian submanifold. We discuss a novel Riemannian geometry for the matrix simplex manifold and show derivation of first- and second-order optimization related ingredients.
\end{abstract}

\section{Introduction} 
Column (or row) stochastic matrices are those where each column (or row) has non-negative entries that sum to $1$. Such matrices are shown to be useful in many machine learning applications \cite{lebanon04a,zhang05a,inokuchi07a,rakotomamonjy08a,sun15a}. The constraint of interest for those matrices is 
\begin{equation}\label{eq:standard_simplex}
x_1 + x_2 + \ldots + x_K = 1, \text{where }x_i \geq 0 \text{ for all } i = \{1,\ldots,K \},
\end{equation}
which is also called the \emph{probability simplex} constraint. 


In this work, we propose to generalize the constraint (\ref{eq:standard_simplex}) to constraints with matrices, i.e., the matrix simplex constraint
\begin{equation}\label{eq:matrix_simplex}
\bX_1 + \bX_2 + \ldots + \bX_K = \bI, 
\end{equation}
where $\bX_i \succeq \bzero $ is a symmetric positive semidefinite matrix of size $n\times n$ for all $i=\{1,2,\ldots, K\}$. Although the constraint (\ref{eq:matrix_simplex}) is a natural generalization of (\ref{eq:standard_simplex}), its study is rather limited \cite{vrehavcek05a,lee08a}. An interesting property of this constraint is its ability to learn mutually orthogonal subspaces with applications in machine learning \cite{giguere2017manifold} and computer vision \cite{kim2009line}.

We discuss a novel Riemannian geometry for the set obtained from the constraint (\ref{eq:matrix_simplex}) with strict positive definiteness of the matrices. Strict positive definiteness of matrices is needed to obtain a differentiable manifold structure. The proposed Riemannian structure allows to handle potentially semidefinite, i.e., rank deficient, matrices gracefully by scaling those elements to the boundary of the manifold. The main aim of this work is to develop optimization-related ingredients that allow Riemannian optimization on this constraint set. The versatile framework of Riemannian optimization enables to perform first-order and second-order optimization of a smooth function. We show the efficacy of the modeling on a toy problem of reconstructing orthogonal subspaces from their noisy estimates.

The manifold related code files are available with the manifold optimization toolbox Manopt \cite{boumal14a}.

\section{The matrix simplex manifold}\label{sec:matrix_simplex}
We define the matrix simplex manifold of interest as 
\begin{equation}\label{eq:manifold}
\begin{array}{lll}
\mathcal{M}_n^K \coloneqq \{(\bX_1 , \bX_2 , \ldots , \bX_K): & \bX_1 + \bX_2 + \ldots + \bX_K = \bI, \\
& \bX_i \in \mathbb{R}^{n\times n}, \text{ and } \\
& \bX_i \succ \bzero \text{ for all }i\in \{1,2,\ldots,K \}  \}.
\end{array}
\end{equation}
It should be noted that the positive semidefiniteness constraint $ \bX_i \succeq \bzero$ is replaced with the positive definiteness constraint $ \bX_i \succ \bzero$ to ensure that the set $\mathcal{M}_n^K  $ is differentiable. Below, we impose a Riemannian structure to the matrix simplex manifold (\ref{eq:manifold}) and discuss ingredients that allow to develop optimization algorithms systematically \cite{absil08a,boumal2020intromanifolds}.

\subsection{Riemannian metric and tangent space projector}
An element $x$ of $\mathcal{M}_n^K$ is numerically represented as the structure $(\bX_1 , \bX_2 , \ldots , \bX_K)$ which is a collection of $K$ symmetric positive definite matrices of size $n\times n$. 

The tangent space of $\mathcal{M}_n^K$ at an element $x$ is the linearization of the manifold, i.e., the constraint~(\ref{eq:matrix_simplex}). Accordingly, the tangent space characterization of $\mathcal{M}_n^K$ at $x$ is
\begin{equation}\label{eq:tangent_space}
\begin{array}{lll}
T_x \mathcal{M}_n^K = \{(\xi_{\bX_1}, \xi_{\bX_2}, \ldots,\xi_{\bX_K}): & \xi_{\bX_1} + \xi_{\bX_2} + \ldots + \xi_{\bX_K} = \bzero  \\
& \xi_{\bX_i} \in \mathbb{R}^{n\times n}, \text{ and } \\
& \xi_{\bX_i}^\top = \xi_{\bX_i} \text{ for all }i\in \{1,2,\ldots,K \}  \}.
\end{array}
\end{equation}

It can be shown that $\mathcal{M}_n^K$ is an embedding submanifold of the $\mathcal{S}_n^K \coloneqq \mathbb{SPD}_n \times \mathbb{SPD}_n \times \ldots \times_K \mathbb{SPD}_n$, which is the Cartesian product of $K$ manifolds of symmetric positive definite matrices of size $n\times n$ \cite{absil08a,boumal2020intromanifolds}. Here $\mathbb{SPD}_n$ denotes the manifold of $n\times n$ symmetric positive definite matrices that has a well-known Riemannian geometry \cite{bhatia09a}. The dimension of the manifold $\mathcal{M}_n^K$ is $(K-1)n(n+1)/2$. 

We endow the manifold with a smooth metric $g_x: T_x \mathcal{M}_n^K\times T_x \mathcal{M}_n^K \rightarrow \mathbb{R}$ (inner product) at every $x \in \mathcal{M}_n^K$ \cite{absil08a}. A natural choice of the metric is based on the well-known bi-invariant metric of $\mathbb{SPD}_n$ \cite{bhatia09a}, i.e.,
\begin{equation}\label{eq:metric}
g_x(\xi_x, \eta_x) \coloneqq 
\sum\limits_i\trace(\bX_i^{-1} \xi_{\bX_i} \bX_i^{-1} \eta_{\bX_i}).
\end{equation}

Once the manifold $\mathcal{M}_n^K$ is endowed with the metric (\ref{eq:metric}), the manifold $\mathcal{M}_n^K$ turns into a Riemannian submanifold of $\mathcal{S}_n^K$. Following \cite[Chapters~3 and 4]{absil08a}, the Riemannian submanifold structure allows the computation of the Riemannian gradient and Hessian of a function (on the manifold) in a straightforward manner from the partial derivatives of the function.

A critical ingredient in those computations is the computation of the linear projection operator of a vector in the ambient space $\mathbb{R}^{n \times n} \times \mathbb{R}^{n \times n} \times \ldots \times_K \mathbb{R}^{n \times n}$ onto the tangent space (\ref{eq:tangent_space}) at an element of $\mathcal{M}_n^K$. In particular, given $z = (\bZ_1,\bZ_2,\ldots,\bZ_K) \in \mathbb{R}^{n \times n} \times \mathbb{R}^{n \times n} \times \ldots \times_K \mathbb{R}^{n \times n}$ in the ambient space, we compute the projection operator $\Pi_x: \mathbb{R}^{n \times n} \times \mathbb{R}^{n \times n} \times \ldots \times_K \mathbb{R}^{n \times n} \rightarrow T_x \mathcal{M}_n^K$, orthogonal with respect to the metric (\ref{eq:metric}), as \cite{mishra16a}
\begin{equation*}
\Pi_x(z) = \argmin\limits_{\xi_x \in T_x \mathcal{M}_n^K} \quad - g_x(z, \xi_x) + \frac{1}{2} g_x(\xi_x, \xi_x),
\end{equation*}
which has the expression
\begin{equation}\label{eq:projection_operator}
\Pi_x(z)  = (\bZ_1 + \bX_1 {\bf{\Lambda}} \bX_1,\  \bZ_2 + \bX_2 {\bf{\Lambda}} \bX_2,  \ \ldots,\  \bZ_K + \bX_K {\bf{\Lambda}} \bX_K),
\end{equation}
where $\bf{\Lambda}$ is the symmetric matrix that is the solution to the linear system
\begin{equation}\label{eq:linear_system}
\sum\limits_i \bX_i {\bf{\Lambda}} \bX_i = - \sum\limits_i \bZ_i.
\end{equation}
It is easy to verify that
\begin{itemize} 
\item $\Pi_x(z)$ belongs to the tangent space $T_x \mathcal{M}_n^K$ and
\item $z - \Pi_x(z)$ and $\Pi_x(z)$ complementary to each other with respect to the chosen metric (\ref{eq:metric}) for all choices of $z$. 
\end{itemize}

\subsection{Retraction operator}
Given a vector in the tangent space, the retraction operator maps it to an element of the manifold \cite[Chapter~4]{absil08a}. Overall, the notion of retraction operation allows to move on the manifold, which is required by any optimization algorithm. 

A natural choice of the retraction operator on the manifold $\mathcal{M}_n^K$ is inspired from the well-known exponential mapping operation on $\mathbb{SPD}_n$, the manifold of positive definite matrices \cite{bhatia09a}. However, this only ensures positive definiteness of the output matrices. To maintain the summation equal to $\bI$ constraint, we additionally \emph{normalize} in a particular fashion. Overall, given a tangent vector $\xi_x \in T_x \mathcal{M}_n^K$, the expression for the retraction operator $R_x: T_x \mathcal{M}_n^K \rightarrow \mathcal{M}_n^K$ is 
\begin{equation}\label{eq:retraction_operation}
\begin{array}{lll}
R_x(\xi_x) \coloneqq &(  \bY_{\rm sum}^{-1/2}\bY_1 \bY_{\rm sum}^{-1/2},\ 
\bY_{\rm sum}^{-1/2}\bY_2 \bY_{\rm sum}^{-1/2}, \ldots, \ \bY_{\rm sum}^{-1/2}\bY_K \bY_{\rm sum}^{-1/2}),
\end{array}
\end{equation}
where $\xi_x = ( \xi_{\bX_1},  \xi_{\bX_2}, \ldots,  \xi_{\bX_K})$, $\bY_i = \bX_{i}({\rm expm}( \bX_{i}^{-1}  \xi_{\bX_i}))$, $\bY_{\rm sum} = \sum_i \bY_i$, and ${\rm expm}(\cdot)$ is the matrix exponential operator.

To show that the operator (\ref{eq:retraction_operation}) is a retraction operator, we need to verify the conditions \cite[Chapter~4]{absil08a}: 
\begin{itemize}
\item the centering condition, i.e., $R_x(0_x) = x$ and
\item the local rigidity condition, i.e., ${\rm D} R_x(0_x) = {\rm id}_{T_x \mathcal{M}_n^K}$, where ${\rm id}_{T_x \mathcal{M}_n^K}$ denotes the identity mapping on $T_x \mathcal{M}_n^K$. 
\end{itemize}
The centering condition for (\ref{eq:retraction_operation}) is straightforward to verify by setting $\xi_x = 0$. To verify the local rigidity condition, we analyze the differential of the retraction operator locally, which is the composition of two steps: the first one is through the matrix exponential and the second is through the normalization by pre and post multiplying with $\bY_{\rm sum}^{-1/2}$. The matrix exponential is locally rigid due to the fact that it defines the well-known exponential mapping on the $\mathbb{SPD}_n$ manifold \cite{absil08a,bhatia09a}. The normalization step (with pre and post multiplying by $\bY_{\rm sum}^{-1/2}$) does not change local rigidness. Hence, the overall composition (\ref{eq:retraction_operation}) satisfies both the centering and local rigidity conditions needed to be a retraction operation.

\subsection{Riemannian gradient and Hessian computations}

As mentioned earlier, a benefit of the Riemannian submanifold structure is that it allows to compute the Riemannian gradient and Hessian of a function in a systematic manner. To that end, we consider a smooth function $f:\mathcal{M}_n^K \rightarrow \mathbb{R}$ on the manifold. We also assume that it is well-defined on $\mathcal{S}_n^K$.

If $\nabla_x f$ is the Euclidean gradient of $f$ at $x \in \mathcal{M}_n^K$, then the Riemannian gradient $\grad_x f$ has the expression
\begin{equation*}
\begin{array}{lll}
\grad_x f & = \Pi_x(\text{gradient on } \mathcal{S}_n^K) \\
& = \Pi_x(\bX_1 {\rm symm}(\nabla_{\bX_1} f)\bX_1, \ \bX_2 {\rm symm}(\nabla_{\bX_2} f)\bX_2,\ \ldots, \ \bX_K {\rm symm}(\nabla_{\bX_K} f)\bX_K ),
\end{array}
\end{equation*}
where $\nabla_{\bX_i}f$ is the partial derivative of $f$ at $x$ with respect to $\bX_i$ and $\Pi_x$ is the tangent space projection operator defined in (\ref{eq:projection_operator}). Here, ${\rm symm}(\cdot)$ extracts the symmetric part of a matrix, i.e., ${\rm symm}({\bf \Delta}) = ( {\bf \Delta} + {\bf \Delta} ^\top)/2$.

The computation of the Riemannian Hessian on the manifold $\mathcal{M}_n^K$ involves the notion of Riemannian connection \cite[Section~5.5]{absil08a}. The Riemannian connection, denoted as  $\nabla_{\xi_{x}}\eta_{x}$, at $x \in \mathcal{M}_n^K$ generalizes the \emph{covariant-derivative} of the tangent vector $\eta_{x} \in T_x \mathcal{M}_n^K$ along the direction of the tangent vector $\xi_{x} \in T_x \mathcal{M}_n^K$ on the manifold $\mathcal{M}_n^K$. Since $\mathcal{M}_n^K$ is a Riemannian submanifold of the manifold $\mathcal{S}_n^K$, the computation of the Riemannian connection enjoys a simple expression in terms of the computations on the symmetric positive definite manifold $\mathbb{SPD}_n$ \cite{bhatia09a}. In particular, the Riemannian connection on $\mathcal{M}_n^K$ is obtained by restricting the connection on $\mathcal{S}_n^K$ to the tangent space $T_x \mathcal{M}_n^K$. The connection on $\mathcal{S}_n^K$ is easy to derive thanks to the well-known Riemannian geometry of $\mathbb{SPD}_n$. Overall, the Riemannian connection expression for $\mathcal{M}_n^K$ is 
\begin{equation} \label{eq:Riemannian_connection}
\begin{array}{lll}
\nabla_{\xi_{x}}\eta_{x} & = \Pi_x (\text{connection on }\mathcal{S}_n^K)\\
                 & = \Pi_x ( {\rm D}\eta_{x}[\xi_{x}] - ({\rm symm}(\xi_{\bX_1}  \bX_1^{-1} \eta_{\bX_1}),  \ldots,  {\rm symm}(\xi_{\bX_K}  \bX_K^{-1} \eta_{\bX_K})   )),
\end{array}
\end{equation}
where ${\rm D} \eta_{x}[\xi_{x}]$ denotes the directional derivative of $\eta_x$ along $\xi_x$. Based on the expression (\ref{eq:Riemannian_connection}), the Riemannian Hessian operation $\hess_x f [\xi_x]$ along a tangent vector $\xi_x \in T_x \mathcal{M}_n^K$ has the expression 
\begin{equation*}\label{eq:Riemannian_Hessian}
\begin{array}{lll}
\hess_x f [\xi_x] & =  \nabla_{\xi_{x}}\grad_x f.
\end{array}
\end{equation*}

\subsection{Computational cost of optimization ingredients}\label{sec:computational_cost}
The expressions shown earlier involve matrix operations that cost $O(n^3 K)$. The solution to the system (\ref{eq:linear_system}) can be obtained iteratively using standard linear equation solvers. The overall cost for the computations is linear in $K$.

\section{Extension to other cases}

\subsection*{The Hermitian case}
The earlier developments easily extend to Hermitian positive definite matrices satisfying the constraint (\ref{eq:matrix_simplex}). The matrix transpose operation is replaced with the conjugate transpose operation \cite{sra15a}. The expressions in Section \ref{sec:matrix_simplex} are modified accordingly.

\subsection*{Decomposition of a symmetric positive definite matrix}
The constraint (\ref{eq:matrix_simplex}) can be generalized to decomposition of a general symmetric positive definite matrix $\bM \succ \bzero$, i.e.,
\begin{equation}\label{eq:matrix_simplex_general}
\bX_1 + \bX_2 + \ldots + \bX_K = \bM. 
\end{equation}
The constraint (\ref{eq:matrix_simplex_general}) can be translated to the simplex constraint (\ref{eq:matrix_simplex}) with a variable change by pre- and post-multiplying with $\bM^{-1/2}$. Consequently, the optimization ingredients developed in Section \ref{sec:matrix_simplex} can be used to perform optimization over the constraint set (\ref{eq:matrix_simplex_general}).

\subsection*{Dealing with large-scale scenarios}
As discussed in Section \ref{sec:computational_cost}, the manifold related operations cost $O(n^3)$. For large values of $n$, the computational cost is prohibitive. A way forward is to learn the matrices $\bX_i$s in a restricted subspace. This is achieved by modeling $\bX_i = \bU \bB_i \bU^\top$ as a rank-$r$ matrix, where the subspace is captured by $\bU \in \mathbb{R}^{n\times r}$ and $\bB_i \succ \bzero$ is a $r$-by-$r$ symmetric positive definite matrix \cite{bonnabel2010riemannian}. It should be noted that $\bU$ is common across all $\bX_i$s. Learning of $\bU$ and $\bB_i$s is then equivalent to optimizing with the constraints $\bU^\top\bU = \bI$ and $\bB_1 + \bB_2 + \ldots + \bB_K = \bI$, which is viewed as a Cartesian product of Stiefel and matrix simplex manifolds. The computations of optimization ingredients in this case cost $O(nr^2 + r^3K)$ which is much less than $O(n^3 K)$ for $r \ll n$.

\section{Toy example: learning mutually orthogonal subspaces}
The constraint (\ref{eq:matrix_simplex}) allows to learn mutually orthogonal subspaces, i.e., $\bX_1$, $\bX_2$, \ldots, and $\bX_K$ are mutually orthogonal. To this end, we consider three matrices $\bA_1$, $\bA_2$, and $\bA_3$ of size $100\times 100$ such that they share mutually orthogonal eigenvectors and have all eigenvalues as either $1$ or $0$. The eigenvalues are added with uniformly random numbers between $0$ and $1$ to generate noisy matrices $\widehat{\bA}_1$,  $\widehat{\bA}_2$, and  $\widehat{\bA}_3$. We solve the problem:
\begin{equation}\label{eq:toy_problem}
\min\limits_{(\bX_1 , \bX_2 , \ldots , \bX_K) \in \mathcal{M}^K_n} - \sum\limits_{k} \trace(\bX_k \widehat{\bA}_k).
\end{equation}
Here, $n = 100$ and $K=3$. To measure effectiveness of learning mutually orthogonal matrices, we define the mutual orthogonality error for $\bX_i$ as $\frac{\trace(\bX_i (\bI - \bX_i))}{\|\bX_i \|_{\rm Fro} \|\bI - \bX_i \|_{\rm Fro}}$. A zero mutual orthogonality error for $\bX_i$ ensures that $\bX_i$ is orthogonal to other matrices.

\begin{figure}
	\begin{center}
		\begin{minipage}[t]{.3\textwidth}
			\begin{center}
				\includegraphics[width=\textwidth]{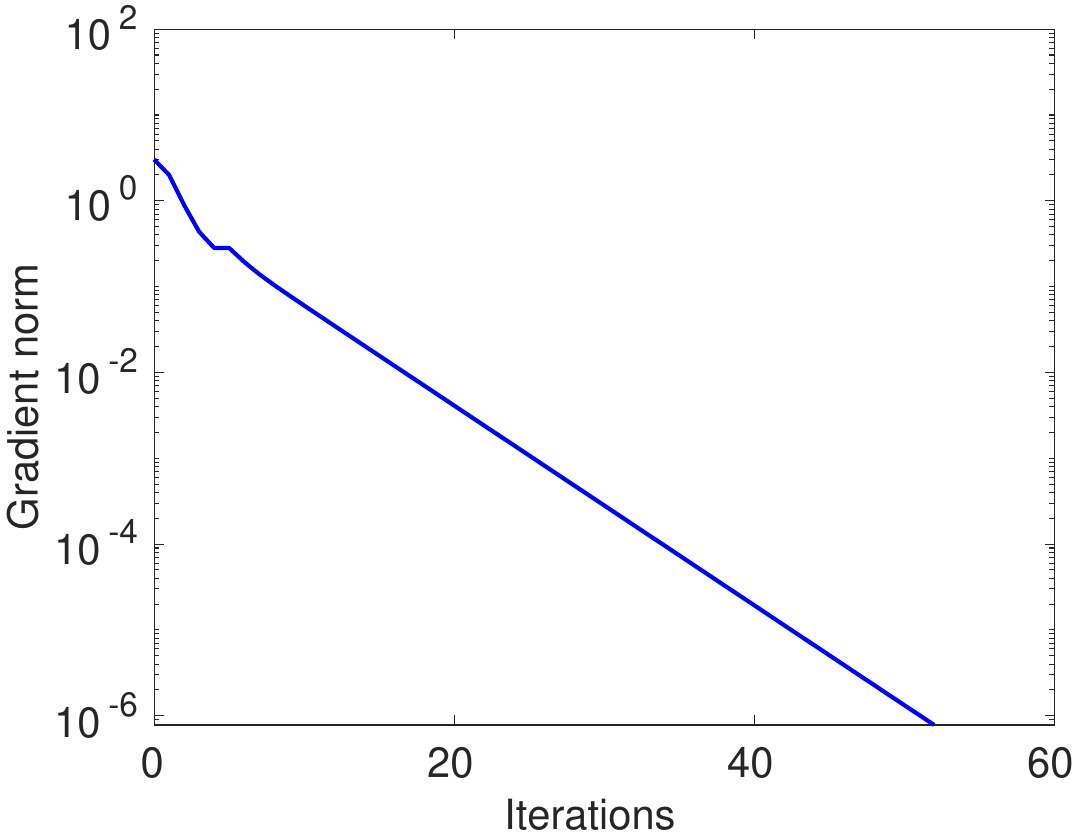}\\	
			\end{center} 
		\end{minipage}
		\begin{minipage}[t]{.3\textwidth}
			\begin{center}
				\includegraphics[width=\textwidth]{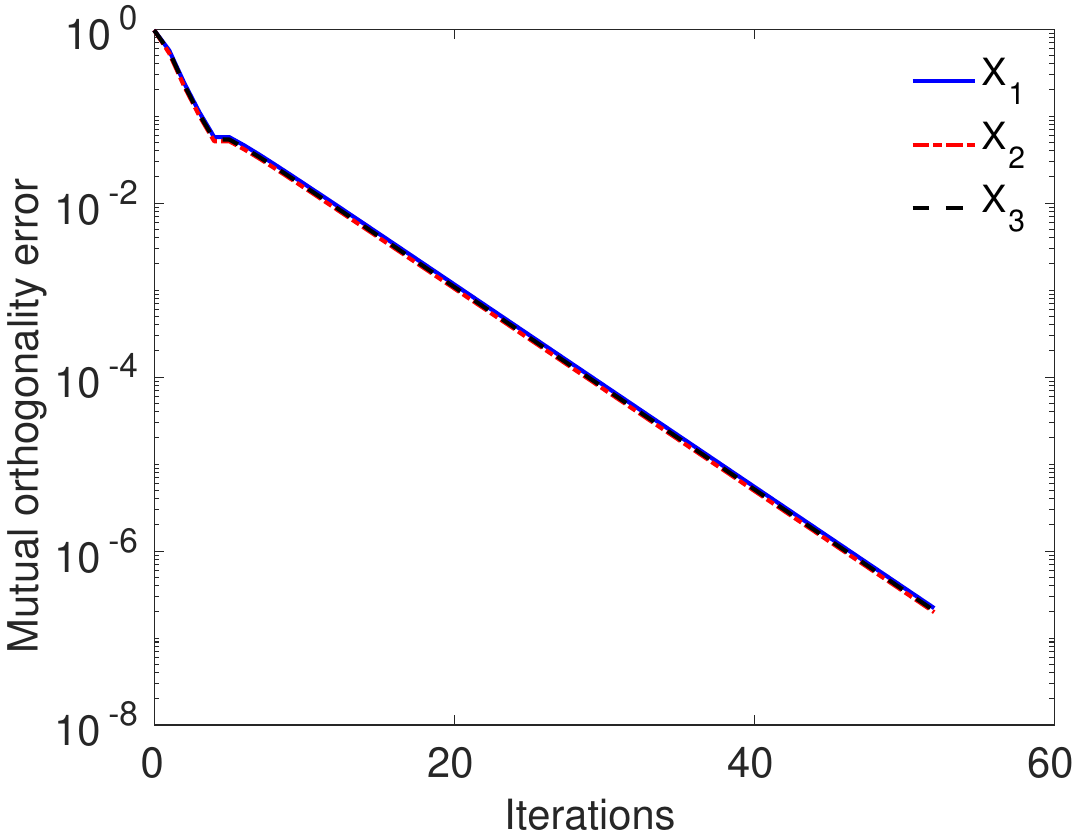}\\		
			\end{center} 
		\end{minipage}
		\begin{minipage}[t]{.3\textwidth}
			\begin{center}
				\includegraphics[width=\textwidth]{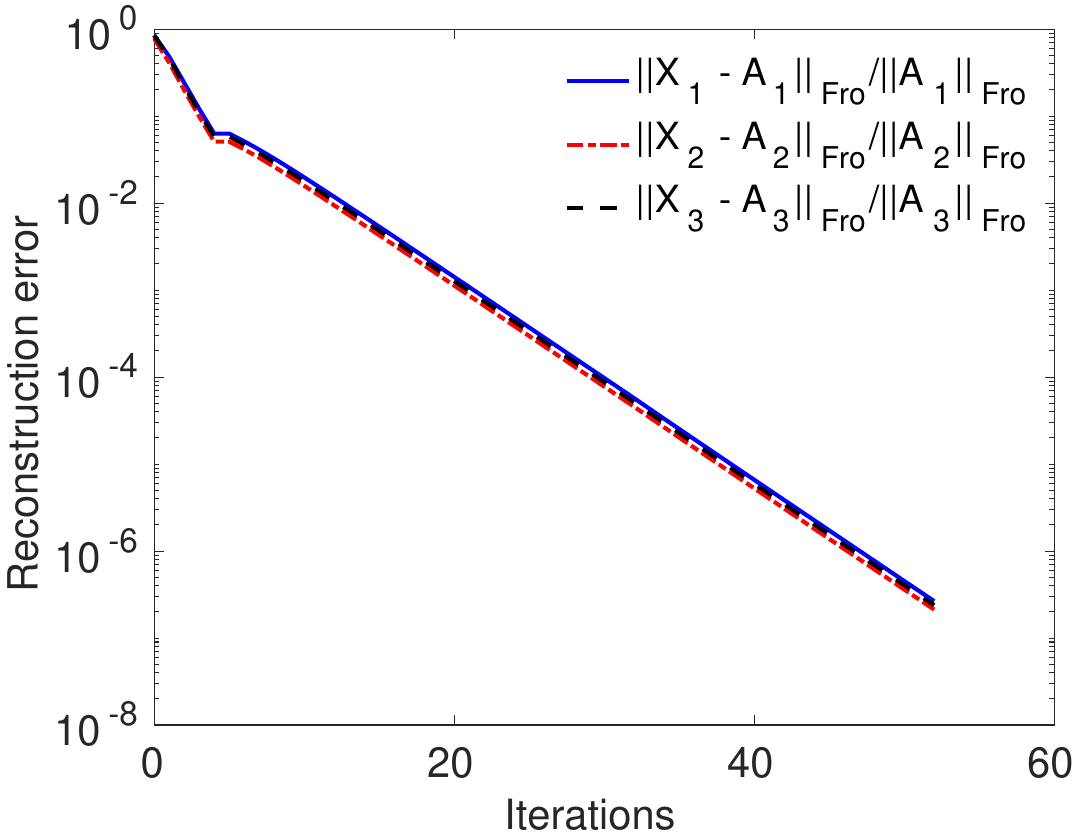}\\			
			\end{center} 
		\end{minipage}
		\caption{The matrix simplex constraint allows learning of mutually orthogonal subspaces. Solving (\ref{eq:toy_problem}) leads to identifying mutually orthogonal matrices while ensuring a low reconstruction error.}
		\label{fig:mutually_orthogonal}
	\end{center}
\end{figure}

We implement a trust-region algorithm for solving (\ref{eq:toy_problem}). Figure \ref{fig:mutually_orthogonal} shows the performance of the algorithm on the problem instance. We observe that our algorithm is able to achieve a high degree of mutual orthogonality between the learned matrices as the error falls below $10^{-6}$. 

In this case, the mutual orthogonality error for $\widehat{\bA}_1$, which is $\frac{\trace(\widehat{\bA}_1 (\widehat{\bA}_2 + \widehat{\bA}_3))}{\|\widehat{\bA}_1 \|_{\rm Fro} \|\widehat{\bA}_2 + \widehat{\bA}_3 \|_{\rm Fro}}$ is $0.68$. Similarly, the  mutual orthogonality errors for $\widehat{\bA}_2$ and $\widehat{\bA}_3$ are $0.66$ and $0.65$, respectively.

\section{Conclusion}
We discussed the matrix simplex manifold as a generalization of the popular probability simplex constraint to symmetric positive semidefinite matrices. The main aim of the work was to understand the geometry of the manifold from an optimization perspective. To this end, the expressions of optimization-related tools were developed. 

\section*{Acknowledgment}

H. Kasai was partially supported by JSPS KAKENHI Grant Numbers JP16K00031 and JP17H01732, and by Support Center for Advanced Telecommunications Technology Research (SCAT).

\bibliography{myreferences}

\end{document}